\documentclass[11pt]{amsart}

\usepackage[a4paper,hmargin=3.5cm,vmargin=4cm]{geometry}
\usepackage{amsfonts,amssymb,amscd,amstext}
\usepackage{graphicx}
\usepackage[dvips]{epsfig}

\usepackage{fancyhdr}
\usepackage{times}

\usepackage{enumerate}
\usepackage{titlesec}
\usepackage{mathrsfs}

\titleformat{\section}
{\filcenter\bfseries\large} {\thesection{.}}{0.2cm}{}
\titleformat{\subsection}[runin]
{\bfseries} {\thesubsection{.}}{0.15cm}{}[.]
\titleformat{\subsubsection}[runin]
{\em}{\thesubsubsection{.}}{0.15cm}{}[.]

\usepackage[up,bf]{caption}

\newtheorem{theorem}{Theorem}[section]

\theoremstyle{definition}

\numberwithin{equation}{section}
\numberwithin{figure}{section}

\def\Tcal{\mathcal{T}}

\def\Cscr{\mathscr{C}}

\def\c{\mathbb{C}}
\def\z{\mathbb{Z}}
\def\d{\mathbb{D}}

\def\r{\mathbb{R}}
\def\n{\mathbb{N}}
\def\t{\mathbb{T}}
\def\z{\mathbb{Z}}
\def\h{\mathbb{H}}

\def\Agot{\mathfrak{A}}
\def\Egot{\mathfrak{E}}

\def\dist{\mathrm{dist}}

\newcommand\wt{\widetilde}

\usepackage{color}

\begin{document}

\fancyhead[LO]{Null holomorphic curves in $\c^3$ and the
conformal Calabi-Yau problem}
\fancyhead[RE]{A.\ Alarc\'on and F.\ Forstneri\v c}
\fancyhead[RO,LE]{\thepage}

\thispagestyle{empty}

\vspace*{1cm}
\begin{center}
{\bf\LARGE Null holomorphic curves in $\c^3$ and applications to the
conformal Calabi-Yau problem}

\vspace*{0.5cm}

{\large\bf Antonio Alarc\'on $\;$ and $\;$ Franc Forstneri\v c}
\end{center}


\vspace*{1cm}

\begin{quote}
{\small
\noindent {\bf Abstract}\hspace*{0.1cm} 
In this paper we survey some recent contributions by the authors \cite{AF1,AF2,AF3} to the theory of null holomorphic curves in the complex Euclidean space $\c^3$, as well as their applications to null holomorphic curves in the special linear group $SL_2(\c)$, minimal surfaces in the Euclidean space $\r^3$, and constant mean curvature one surfaces (Bryant surfaces) in the hyperbolic space $\h^3$. The paper is an expanded version of the lecture given by the second named author at the Abel Symposium 2013 in Trondheim.

\vspace*{0.1cm}

\noindent{\bf Keywords}\hspace*{0.1cm} Riemann surfaces, complex curves, null holomorphic curves, minimal surfaces, Bryant surfaces,  complete immersions, proper immersions.

\vspace*{0.1cm}

\noindent{\bf MSC (2010):}\hspace*{0.1cm} 53C42; 32H02, 53A10, 32B15.}
\end{quote}


\section{On null curves in $\c^3$ and minimal surfaces in $\r^3$}
\label{sec:C3}
An open connected Riemann surface $M$ is said to be a {\em bordered Riemann surface} if it
is the interior of a compact one dimensional complex manifold $\overline M$ with smooth boundary
$bM\neq\emptyset$ consisting of finitely many closed Jordan curves. The closure $\overline M = M\cup bM$ of such $M$ is a {\em compact bordered Riemann surface}. By classical results every compact bordered Riemann surface is conformally equivalent to a smoothly bounded domain in an open (or compact) Riemann surface $R$ by a map smoothly extending to the boundary.

Let $M$ be an open Riemann surface. A {\em holomorphic} immersion $F=(F_1,F_2,F_3)\colon M\to\c^3$ is said to be a {\em null curve} if it is directed by the conical quadric subvariety
\begin{equation}\label{eq:null}
		\Agot=\{(z_1,z_2,z_3)\in\c^3\colon z_1^2+z_2^2+z_3^2=0\}, 
\end{equation}
in the sense that the derivative $F'=(F_1',F_2',F_3')\colon M\to\c^3$ with respect to any local holomorphic coordinate on $M$ takes values in $\Agot\setminus\{0\}$. If $M$ is a bordered Riemann surface, then the same definition applies to smooth maps $\overline M\to\c^3$ which are holomorphic in $M$.

The real and the imaginary part of a null curve $M\to\c^3$ are {\em conformal} (angle preserving) {\em minimal immersions} $M\to\r^3$; that is, having mean curvature identically zero. Conversely, every conformal minimal immersion $M\to\r^3$ is locally (on any simply-connected domain in $M$) the real part of a null curve; this fails on non-simply connected Riemann surfaces due to the periods of the harmonic conjugate. This connection enables the use of complex analytic tools in minimal surface theory, a major topic of differential geometry since the times of Euler and Lagrange. One of the quintessential examples of this statement is the use of Runge's approximation theorem for holomorphic functions on open Riemann surfaces in the study of {\em the Calabi-Yau problem for surfaces}.

In 1965, Calabi \cite{Calabi} conjectured the nonexistence of {\em complete} minimal surfaces in $\r^3$ with bounded projection into a straight line; this would imply in particular the nonexistence of bounded complete minimal surfaces in $\r^3$. Recall that an immersion $F\colon M\to\r^n$, $n\in\n$, is said to be complete if the pullback $F^*ds^2$ by $F$ of the Riemannian metric on $\r^n$ is a complete metric on $M$. Calabi's conjecture turned out to be false by the groundbreaking counterexample by Jorge and Xavier \cite{JX} in 1980 who constructed a complete conformal minimal immersion $F=(F_1,F_2,F_3)\colon\d\to\r^3$ of the disc $\d=\{z=x+\imath y\in\c\colon |z|<1\}$ with $F_3(z)=\Re(z)=x$. (Here $\imath=\sqrt{-1}$.) Their method consists of applying the classical Runge theorem on a labyrinth of compact sets in $\d$ in order to construct a suitable harmonic function $(F_1,F_2)\colon\d\to\r^2$. Refinements of Jorge and Xavier's technique have given rise to a number of examples of complete minimal surfaces in $\r^3$ with a bounded coordinate function (see \cite{RT,Lopez1,Lopez2}); in particular it was recently proven by Alarc\'on, Fern\'andez, and L\'opez \cite{Afer,AFL1,AFL2} that every open Riemann surface carrying non-constant bounded harmonic functions is the underlying conformal structure of such a surface.

As pointed out by S.-T.\ Yau in his 1982 problem list \cite[Problem 91]{Yau1}, the question of whether there exist complete bounded minimal surfaces in $\r^3$ (which became known in the literature as {\em the Calabi-Yau problem}) remained open. (See also Yau's {\em Millenium Lecture 2000} \cite{Yau2}.) It was Nadirashvili \cite{Nad} who in 1996 answered this question in the affirmative by constructing a conformal complete bounded minimal immersion of the disc $\d$ into $\r^3$. Nadirashvili's method relies on a recursive use of Runge's approximation theorem on labyrinths of compact subsets of $\d$, and it has been the seed of several construction techniques, leading to a variety of examples. In particular, in 2012 Ferrer, Mart\'in, and Meeks \cite{FMM} constructed complete properly immersed minimal surfaces with  {\em arbitrary topology} in any given either convex or bounded and smooth domain of $\r^3$; see also Alarc\'on, Ferrer, and Mart\'in \cite{AFM} for the case of finite topology. 

A question that appeared early in the history of the Calabi-Yau problem was whether there exist complete bounded minimal surfaces in $\r^3$ whose conjugate surfaces also exist and are bounded; that is, whether there exist complete bounded null curves in $\c^3$. The first such examples were provided only very  recently by Alarc\'on and L\'opez \cite{AL} who constructed complete null curves with arbitrary topology properly immersed in any given convex domain of $\c^3$; this answers a question by Mart\'in, Umehara, and Yamada \cite[Problem 1]{MUY1}. Their method, which is different from Nadirashvili's one, relies on a Runge-Mergelyan type theorem for null curves in $\c^3$ \cite{AL1}, a new and powerful tool that gave rise to a number of constructions of both minimal surfaces in $\r^3$ and null curves in $\c^3$ (see \cite{AL1,AFL1,AL,AL-Israel,AL-conjugate}). Very recently, Ferrer, Mart\'in, Umehara, and Yamada \cite{FMUY} showed that Nadirashvili's method can be adapted to null curves, giving an alternative proof of the existence of complete bounded null discs in $\c^3$. 

In the opposite direction, Colding and Minicozzi \cite{CM} proved in 2005 that every complete {\em embedded} minimal surface of finite topology in $\r^3$ is proper in $\r^3$, hence unbounded. This result was extended by Meeks, P\'erez, and Ros \cite{MPR} to surfaces with finite genus and countably many ends. It follows that the original Calabi's conjecture is true for embedded surfaces of finite topology. Although being embedded is a strong constraint for a complete minimal surface in $\r^3$, it is no constraint for null curves in $\c^3$ as the following result shows.

\begin{theorem}[\text{\cite[Corollary 6.2]{AF2}}]\label{th:topologyembedded}
There exist complete null curves with arbitrary topology properly embedded in any given convex domain of $\c^3$.
\end{theorem}

This answers a question by Mart\'in, Umehara, and Yamada \cite[Problem 2]{MUY1}. 

The key to the proof is that {\em the general position of null curves in $\c^3$ is embedded}. In fact, Theorem \ref{th:topologyembedded} easily follows from the existence of complete properly immersed null curves in convex domains of $\c^3$ \cite{AL} and the following desingularization result from \cite{AF2}.

\begin{theorem}[\text{\cite[Theorem 2.4]{AF2}}]\label{th:generalposition}
Let $M$ be a bordered Riemann surface. Then every  null immersion $\overline M\to\c^3$ can be approximated in the $\Cscr^1$-topology on $\overline M$ by null embeddings $\overline M\hookrightarrow\c^3$. Similarly, if $M$ is any open Riemann surface then any null immersion $M\to \c^3$ can be approximated uniformly on compacts by null embeddings $M \to\c^3$.
\end{theorem}

Since complex submanifolds of complex Euclidean spaces are area minimizing \cite{Federer}, the Calabi-Yau problem is closely related to a question, posed by Yang \cite{Yang1,Yang2} in 1977, whether there exist complete bounded complex submanifolds of $\c^n$ for $n>1$. The first result in this subject was obtained by Jones \cite{Jones} who constructed  holomorphic immersions $\d\to\c^2$ and embeddings $\d\hookrightarrow\c^3$ with bounded image and complete induced metric. His method is based on the BMO duality theorem. In 2009, Mart\'in, Umehara, and Yamada \cite{MUY2} extended Jones' result to complete bounded complex curves in $\c^2$ with arbitrary finite genus. Finally, Alarc\'on and L\'opez \cite{AL} constructed complete complex curves with arbitrary topology properly immersed in any given convex domain of $\c^2$, as well as the first example of a complete bounded {\em embedded} complex curve in $\c^2$ \cite{AL-embedded}.

The use of Runge's theorem in the recursive procedure of Nadirashvili's method for constructing complete bounded minimal surfaces does not enable one to control the placement in $\r^3$ of the entire surface at each step. Therefore one is forced to cut away some small pieces of the surface in order to keep it suitably bounded, so it is impossible to control the conformal structure on the surface when applying this technique to non-simply connected Riemann surfaces (the simply-connected ones are of course conformally equivalent to the disc $\d$). This phenomenon has been present in every construction of complete bounded minimal surfaces in $\r^3$ and null curves in $\c^3$ (see \cite{AL} and references therein), and also in every construction of complete bounded complex curves in $\c^2$ with non-trivial topology up to this point (cf. \cite{MUY2,AL,AL-embedded}). 

This constraint has recently been overcome by the authors \cite{AF1} in the case of bordered Riemann surfaces in $\c^n$, $n\geq 2$.

\begin{theorem}[\text{\cite[Theorem 1]{AF1}}]\label{th:CY-C2}
Every bordered Riemann surface carries a complete proper holomorphic immersion to the unit ball of $\c^n$,  $n\geq 2$, which can be chosen an embedding for $n\geq 3$.
\end{theorem}

The construction in the proof of Theorem \ref{th:CY-C2} is inspired by that of Alarc\'on and L\'opez \cite{AL}, but it requires additional complex analytic tools. The key point is to replace Runge's theorem by approximate solutions of certain {\em Riemann-Hilbert boundary value problems} in $\c^n$ (see \S\ref{sec:RH} below for an exposition of this subject). This gives sufficient control of the placement of the whole curve in the space to avoid shrinking, thereby enabling one to control its conformal structure. Another important tool is the method of Forstneri\v c and Wold \cite{FW1} for {\em exposing boundary points} of a complex curve. This method, together with a local version of the Mergelyan theorem, enables one to incorporate suitable arcs to a compact bordered Riemann surface in $\c^n$ without modifying its conformal structure.

However, the case of minimal surfaces in $\r^3$ and null curves in $\c^3$ is still much harder and requires more refined complex analytic tools. The following version of the {\em Riemann-Hilbert problem for null curves} has been obtained recently in \cite{AF3}.

\begin{theorem}[\text{\cite[Theorem 3.4]{AF3}}]\label{th:RH}
Let $M$ be a bordered Riemann surface, and let $I$ be a compact subarc of $bM$ which is not a connected component of $bM$. Choose a small annular neighborhood $A\subset \overline M$ of $bM$ and a smooth retraction $\rho\colon A\to bM$. Let $F\colon \overline M \to\c^3$ be a null holomorphic immersion, let $\vartheta\in \Agot\setminus\{0\}$ be a null vector, let $\mu \colon bM \to \r_+$ be a continuous function supported on $I$, and consider the continuous map  
\[
	\varkappa\colon bM \times\overline{\d}\to\c^3,\qquad 
	\varkappa(x,\xi)=F(x) + \mu(x)\,\xi \,\vartheta.
\]
Then for any number $\epsilon>0$ there exist an arbitrarily small open neighborhood $\Omega$ of $I$ in $A$ and a null holomorphic immersion $G\colon \overline M\to\c^3$  satisfying the following properties:
\begin{itemize}
\item[\rm (a)] $\dist(G(x),\varkappa(x,b\d))<\epsilon$ for all $x\in bM$.
\item[\rm (b)] $\dist(G(x),\varkappa(\rho(x),\overline{\d}))<\epsilon$ for all $x\in \Omega$.
\item[\rm (c)] $G$ is $\epsilon$-close to $F$ in the $\Cscr^1$ topology on $\overline M\setminus \Omega$.
\end{itemize}
\end{theorem}

The authors gave a direct proof by explicit calculation in the special case when $M$ is the disc $\d$, using the so called spinor representation of the null quadric ({\ref{eq:null}). This can be used locally on small discs abutting the boundary of $M$. The proof of the general case depends on the technique of {\em gluing holomorphic sprays} (which amounts to a nonlinear version of the $\overline\partial$-problem in complex analysis, see \cite[Chapter 5]{F2011}) applied to the derivatives of null curves. To use the method of gluing sprays in the present setting, one must control the periods of some maps in the amalgamated spray in order to get well defined null curves by integration; in addition, delicate estimates are needed to ensure that the resulting null curves have the desired properties. A more precise outline of the proof is given in \S\ref{sec:RH2}.

The following result, whose proof uses the techniques described above, is the authors' main contribution to the Calabi-Yau problem.

\begin{theorem}[\text{\cite[Theorem 1.1]{AF3}}]\label{th:CY-C3}
Every bordered Riemann surface carries a complete proper null embedding into the unit ball of $\c^3$.
\end{theorem}

If $F\colon M\to\c^3$ is a null curve, then the Riemannian metric $F^*ds^2$ induced by the Euclidean metric of $\c^3$ via $F$ is twice the one 
induced by the Euclidean metric of $\r^3$ via the real part $\Re F$ (cf.\ \cite{Osserman}). Indeed, write $F=(F^1,F^2,F^3)\colon M\to\c^3$ and let $\zeta=x+\imath y$ be a local holomorphic coordinate on $M$. Then
\begin{eqnarray*}
0 & = & \sum_{j=1}^3 (F^j_\zeta)^2=\sum_{j=1}^3 (F^j_x)^2=\sum_{j=1}^3 \big((\Re F^j)_x + \imath (\Im F^j)_x\big)^2\\
 & = & \sum_{j=1}^3 \big(((\Re F^j)_x)^2 - ((\Im F^j)_x)^2\big) + 2\imath \sum_{j=1}^3 (\Re F^j)_x (\Im F^j)_x.
\end{eqnarray*}
Since $(\Im F)_y=-(\Re F)_x$ by the Cauchy-Riemann equations, the above is equivalent to $|(\Re F)_x|=|(\Im F)_y|$ and $\langle(\Re F)_x,(\Im F)_y\rangle=0$. It follows that the minimal immersion $\Re F\colon M\to\r^3$ is conformal, harmonic, and
\[
F^*ds^2=|F_x|^2(dx^2+dy^2)=2|(\Re F)_x|^2(dx^2+dy^2)=2(\Re F)^*ds^2.
\]

In particular, the real part of a complete null curve in $\c^3$ is a complete conformal minimal immersion in $\r^3$. In view of Theorem \ref{th:CY-C3}, we obtain the following result on the so-called {\em conformal Calabi-Yau problem} for surfaces.
\begin{theorem}[\text{\cite[Corollary 1.2]{AF3}}]\label{th:CY-R3}
Every bordered Riemann surface $M$ carries a conformal complete minimal immersion $M\to\r^3$ with bounded image.
\end{theorem}

We emphasize that, in Theorems \ref{th:CY-C3} and \ref{th:CY-R3}, {\em the conformal structure on the source Riemann surface is not changed}. 

Notice that we do not control the asymptotic behavior of the surfaces in Theorem \ref{th:CY-R3}; in particular we do not know whether there exist conformal complete {\em proper} minimal immersions from any bordered Riemann surface into a ball of $\r^3$.

At this point we wish to draw a certain analogy with the old problem whether every open Riemann surface admits a proper holomorphic embedding into $\c^2$ (see Bell and Narasimhan \cite[Conjecture 3.7, p.\ 20]{BN}). It is classical that every open Riemann surface properly holomorphically embeds in $\c^3$ and immerses in $\c^2$ \cite[\S 8.2]{F2011}. By introducing the technique of exposing boundary points alluded to above, combined with the Anders\'en-Lempert theory of holomorphic automorphisms of $\c^n$ for $n>1$, Forstneri\v c and Wold proved in 2009 that, if a compact bordered Riemann surface $\overline M$ admits a (nonproper) holomorphic embedding in $\c^2$  then its interior $M$ admits a {\em proper} holomorphic embedding in $\c^2$ \cite{FW1}. Further applications of their technique can be found in \cite{Majcen} and \cite{FW2}. For example, {\em every circular domain in the Riemann sphere admits a proper holomorphic embedding in $\c^2$} \cite{FW2}. (The case of finitely connected plane domains was established by Globevnik and Stens\o nes in 1995 \cite{GS}.) The main novelty in \cite{FW1,FW2} is that, unlike in the earlier constructions, no cutting of the surface is needed and hence the conformal structure is preserved.  It is considerably easier to show that every open oriented real surface $M$ admits a complex structure $J$ such that the open Riemann surface $(M,J)$ admits a proper holomorphic embedding into $\c^2$ (Alarc\'on and L\'opez \cite{AL3}; for the case of finite topology see \cite{CF}).

Another good example of how Runge's theorem has been exploited in minimal surface theory is the construction of proper minimal surfaces in $\r^3$ with hyperbolic conformal structure. (An open Riemann surface is said to be {\em hyperbolic} if it carries negative non-constant subharmonic functions; otherwise it is called {\em parabolic}.) An old conjecture of Sullivan \cite{Meeks} asserted that every properly immersed minimal surface in $\r^3$ with finite topology must have parabolic conformal structure. The first counterexample was given by Morales \cite{Morales} in who in 2003 constructed a proper {\em conformal} minimal immersion $\d\to\r^3$. 

In the same line, Schoen and Yau \cite{SY} asked in 1985 whether a minimal surface in $\r^3$ properly projecting into a plane must be parabolic. The question was answered by Alarc\'on and L\'opez \cite{AL1} who showed that in fact every open Riemann surface $M$ carries a conformal minimal immersion $X=(X_1,X_2,X_3)\colon M\to\r^3$ such that $(X_1,X_2)\colon M\to\r^2$ is a proper map, and a null curve $F=(F_1,F_2,F_3)\colon M\to\c^3$ such that $(F_1,F_2)\colon M\to\c^2$ is proper. (See also \cite{AL-conjugate}.) Taking into account Theorem \ref{th:generalposition} we obtained the following extension of the previous results.

\begin{theorem}[\text{\cite[Theorem 8.1]{AF2}}]\label{th:properC3}
Every open Riemann surface $M$ carries a proper null {\em embedding} $F=(F_1,F_2,F_3)\colon M\hookrightarrow \c^3$ such that $(F_1,F_2)\colon M\to\c^2$ is a proper map. 
\end{theorem}

The constructions in \cite{AL1} and \cite{AF2} involved particular versions of Runge's theorem for minimal surfaces in $\r^3$ and null curves in $\c^3$ that do not give any control on the third coordinate. However, the newly developed complex analytic methods involved in the proof of Theorem \ref{th:CY-C3} above enable us to overcome this constraint in the case of bordered Riemann surfaces as null curves in $\c^3$. The following is our main result in this line.

\begin{theorem}[\text{\cite[Theorem 1.4]{AF3}}]\label{th:F3}
Every bordered Riemann surface $M$ carries a null holomorphic embedding $F=(F_1,F_2,F_3)\colon M\hookrightarrow \c^3$ such that $(F_1,F_2)\colon M\to\c^2$ is a proper map and the function $F_3\colon M\to\c$ is bounded on $M$.
\end{theorem}

Joining together the methods in the proof of the above result and the Mergelyan theorem for null curves \cite{AL1} (see also \cite{AF2}), we also get the following result.

\begin{theorem}[\text{\cite[Theorem 1.8]{AF3}}]\label{th:F3'}
Every orientable noncompact smooth real surface $M$ without boundary admits a complex structure $J$ such that the Riemann surface $(M,J)$ carries a proper holomorphic null embedding $(F_1, F_2, F_3)\colon (M, J) \to \c^3$ such that $F_3$ is a bounded function on $M$.
\end{theorem}

Recall that the only properly immersed minimal surfaces in $\r^3$ with a bounded coordinate function are planes by the Half-Space theorem of Hoffman and Meeks \cite{HM}. Therefore, if $F=(F_1,F_2,F_3)\colon M\hookrightarrow \c^3$ is a proper null curve as those in Theorem \ref{th:F3}, then $\Re(e^{\imath t} F)\colon M\to\r^3$ is a conformal complete {\em non-proper} minimal immersion for all $t\in\r$.


\section{On null curves in $SL_2(\c)$ and Bryant surfaces in $\h^3$}
\label{sec:SL2C}

A {\em holomorphic} immersion $F\colon M\to SL_2(\c)$ from an open Riemann surface into the special linear group 
\[
	SL_2(\c)= \left \{ z= \left(\begin{matrix} z_{11} & z_{12} 
		\cr z_{21} & z_{22} \end{matrix} \right) \in\c^4 
		\colon \det z= z_{11}z_{22}-z_{12}z_{21}=1\right \}
\]
is said to be a {\em null curve} if it is directed by the quadric variety
\begin{equation}\label{eq:nullSL2C}
\Egot=	\left \{ z= \left(\begin{matrix} z_{11} & z_{12} 
	\cr z_{21} & z_{22} \end{matrix} \right) \colon \det z= z_{11}z_{22}-z_{12}z_{21}=0\right \} 
	\subset \c^4,
\end{equation}
meaning that the derivative $F'\colon M\to\c^4$ with respect to any local holomorphic coordinate on $M$ belongs to $\mathfrak{E}\setminus \{0\}$. If $M$ is a bordered Riemann surface, then the same definition applies for smooth maps $\overline M\to SL_2(\c)$ being holomorphic in $M$.

In 1987, Bryant \cite{Br} discovered that, if $F\colon M\to SL_2(\c)$ is a null curve, then $F\cdot \overline F^T$ takes values in the hyperbolic space $\h^3=SL_2(\c)/SU(2)$ and is a {\em conformal} immersion with constant mean curvature one; conversely, every simply-connected {\em Bryant surface} (i.e., with constant mean curvature one in $\h^3$) is the projection of a null curve in $SL_2(\c)$. As is the case of minimal surfaces in $\r^3$, the above connection enables the use of complex analytic tools in Bryant surface theory which made this subject a fashionable research topic in the last decade (see e.g.\ \cite{UY,CHR,Ro} for the background). Moreover, the {\em Lawson correspondence} \cite{La} implies that every simply connected Bryant surface is isometric to a minimal surface in $\r^3$ and vice versa; hence problems on minimal surface theory are automatically natural also for Bryant surfaces.

Some of the results described in the previous section hold not only for null curves in $\c^3$, but also for immersions with derivative in an arbitrary conical subvariety $A$ of $\c^n$ ($n\geq 3$) which is smooth away from the origin; i.e., $A$-{\em immersions}. 

\begin{theorem}[\text{\cite[Theorem 2.5 and Corollary 2.7]{AF2}}]\label{th:A}
Let $M$ be an open Riemann surface, and let $A$ be a closed conical subvariety of $\c^n$, $n\geq 3$, which is not contained in any hyperplane and such that $A\setminus\{0\}$ is a smooth Oka manifold. Then:
\begin{itemize}
\item (Desingularization theorem) Every $A$-immersion $M\to\c^n$ can be approximated uniformly on compacts by $A$-embeddings $M \to \c^n$.
\item (Runge theorem) Every $A$-immersion $U\to \c^n$ on an open neighborhood $U$ of a compact Runge set $K\subset M$ can be approximated uniformly on $K$ by $A$-immersions $M\to\c^n$.
\end{itemize}
\end{theorem}
Recall that a complex manifold $Y$ is said to be an {\em Oka manifold} if every holomorphic map from an open neighborhood of a compact convex set $K\subset \c^N$ to $Y$ can be approximated, uniformly on $K$, by entire maps $\c^N\to Y$; see \cite{F2011} for a reference on Oka theory.

Although the variety $\Egot$ \eqref{eq:nullSL2C} controlling null curves in $SL_2(\c)$ meets the requirements of Theorem \ref{th:A}, the results do not apply directly since the $\Egot$-immersions $M\to\c^4$ furnished by the theorem need not lie in $SL_2(\c)$. In order to force the image to lie in $SL_2(\c)$, one must add another equation expressing the condition that the tangent vector to the curve is also tangent to $SL_2(\c)$ (as a submanifold of $\c^4$). Unfortunately the resulting system of equations is no longer autonomous (the second equation depends on the point in space), and hence the methods of \cite{AF2} do not apply. 

However, Mart\'in, Umehara, and Yamada \cite{MUY1} discovered in 2009 that the biholomorphic map $\Tcal \colon \c^3\setminus\{z_3=0\} \to SL_2(\c)\setminus\{z_{11}=0\}$, given by
\begin{equation}\label{eq:T}
	\Tcal(z_1,z_2,z_3) = \frac{1}{z_3}\left(\begin{matrix} 1 & z_1+\imath z_2 \cr 
	z_1-\imath z_2 & 	z_1^2+z_2^2+z_3^2
\end{matrix} \right),
\end{equation}
carries null curves into null curves. This transformation allows us to obtain a succulent list of corollaries to the results in the previous section.

In view of Theorem \ref{th:generalposition} we get the following result concerning the general position of null curves in $SL_2(\c)\setminus \{z_{11} = 0\}$.
\begin{theorem}[\text{\cite[Corollary 2.8]{AF2}}]
Let $M$ be a bordered Riemann surface. Every immersed null curve $\overline M\to SL_2(\c)\setminus \{z_{11} = 0\}$ can be approximated in the $\Cscr^1$ topology on $\overline M$ by embedded null curves $\overline M \to SL_2(\c)\setminus \{z_{11} = 0\}$.
\end{theorem}

Observe that for any constant $c>0$ the biholomorphism $\Tcal$ \eqref{eq:T} maps complete bounded null curves in $\c^3\setminus \{|z_3|>c\}$ into complete bounded null curves in $SL_2(\c)$, which in turn project to complete bounded Bryant surfaces in $\h^3$ \cite{MUY1}. From Theorems \ref{th:topologyembedded} and \ref{th:CY-C3}, we get the following results regarding Calabi-Yau type questions.
\begin{theorem}[\text{\cite{AF2}, \cite[Corollary 1.9]{AF3}}]
\item
\begin{itemize}
\item There exist complete bounded embedded null curves in $SL_2(\c)$ and immersed Bryant surfaces in $\h^3$ with arbitrary topology.
\item Every bordered Riemann surface $M$ admits a complete null holomorphic
embedding $M\to SL_2(\c)$ with bounded image, and it is conformally equivalent to a complete
bounded immersed Bryant surface in $\h^3$.
\end{itemize}
\end{theorem}

The latter part of the former item in the above theorem was already proven in \cite{AL1} where also complete bounded {\em immersed} null curves in $SL_2(\c)$ with arbitrary topology were given. Complete bounded immersed simply connected null holomorphic curves in $SL_2(\c)$, hence complete
bounded simply-connected Bryant surfaces in $\h^3$, were provided in \cite{FMUY,MUY1}.

Finally, observe that applying $\Tcal$ to a proper null curve $F = (F_1, F_2, F_3)\colon M \to \c^3$ such that $0<c_1 < |F_3| < c_2$ on $M$ one gets a proper null curve in $SL_2(\c)$, which in turn projects to a proper Bryant surface in $\h^3$. Therefore, Theorems \ref{th:F3} and \ref{th:F3'} provide the following examples of proper null curves in $SL_2(\c)$ and Bryant surfaces in $\h^3$.

\begin{theorem}[\text{\cite[Corollaries 1.5 and 1.6 and Theorem 1.8]{AF3}}]\label{th:SL}
\item
\begin{itemize}
\item There exist properly embedded null curves in $SL_2(\c)$ and properly immersed Bryant surfaces in $\h^3$ with arbitrary topology.
\item Every bordered Riemann surface $M$ admits a proper null holomorphic
embedding $M\to SL_2(\c)$, and it is conformally equivalent to a properly
immersed Bryant surface in $\h^3$.
\end{itemize}
\end{theorem}

Connecting to Sullivan's conjecture for minimal surfaces \cite{Meeks}, the ones in the latter item of Theorem \ref{th:SL} are the first examples of proper null curves in $SL_2(\c)$ and proper Bryant surfaces in $\h^3$ with finite topology and hyperbolic conformal structure.


\section{The Riemann-Hilbert problem and proper holomorphic maps of bordered Riemann surfaces}
\label{sec:RH}
In this and the following section we explain how a certain version of the classical {\em Riemann-Hilbert boundary value problem} is used in the construction of holomorphic curves satisfying various additional properties (for example, proper and/or complete and bounded).

The linear {\em Riemann boundary value problem} on the disc $\d=\{z\in\c:|z|<1\}$ asks for holomorphic functions $f(z)=u(z)+\imath v(z)$ on $\d$, continuous on the closed disc $\overline \d$ and satisfying the following condition on the circle $\t=b\,\d=\{|z|=1\}$:
\begin{equation} \label{eq:R}
	a(z)\, u(z) - b(z)\,v(z) = c(z), \qquad z\in \t,
\end{equation}
where $a$, $b$, and $c$ are given real-valued continuous functions on $\t$. This classical problem considered by Riemann in his dissertation \cite{Riemann}, and its extension to non-simply connected domains and bordered Riemann surfaces, is one of the main boundary value problems of analytic function theory. Geometrically speaking, (\ref{eq:R}) demands that the boundary value of $f$ at any point $z\in \t$ lies in a certain affine real line $l_z\subset \c$ depending on $z$.  

Writing $f(z)=f_+(z)$ and introducing the conjugate function on $\c\setminus\d$ by
\[
	f_-(z) = \overline{f_+\left(\bar{z}^{-1}\right)},\qquad |z|\ge 1,
\]
we have $f_-(z) = \overline{f_+(z)}$ on the circle $\t$, and (\ref{eq:R}) can be written as 
\begin{equation}\label{eq:RH}
    \alpha(z) f_+(z) + \beta(z) f_-(z) = c(z), \qquad |z|=1,
\end{equation}
where $\alpha(z)=\left(a(z)+\imath b(z)\right)/2$ and $\beta(z)=\left(a(z)-\imath b(z)\right)/2$. Hilbert's generalization \cite{Hilbert} asks for functions $f_+$ and $f_-$, holo\-morphic on $\d$ and $\c\setminus \overline\d$, respectively, continuous up to the circle and satisfying (\ref{eq:RH}) where $\alpha$, $\beta$, and $c$ are arbitrary complex-valued functions on $\t$. One may consider the same problem for vector-valued or matrix-valued holomorphic functions. For example, the {\em Birkhoff factorization problem} \cite{Birkhoff} amounts to finding solutions of the equation $f_+(z)=G(z)f_-(z)$ on $|z|=1$, where the matrix-valued functions $f_\pm$ are holomorphic inside and outside of the disc $\d$, respectively.

Riemann's boundary value problem was motivated by the problem of finding a linear differential equation of Fuchsian type with given singular points and a given monodromy group; this is {\em Hilbert's 21st problem} on the famous list of 23 problems from his 1900 ICM lecture. In 1905, Hilbert  \cite{Hilbert} obtained some progress on the general Riemann-Hilbert problem (\ref{eq:RH}) by reducing it to an integral equation. The homogeneous vector-valued Riemann-Hilbert problem was solved in 1908 by Plemelj \cite{Plemelj, Plemelj2} by using his {\em jump formula} for Cauchy integrals. This also gave a solution of Hilbert's 21st problem, except for a gap in some cases (see Bolibrukh \cite{Boli} and Kostov \cite{Kostov}). In 1909 Birkhoff \cite{Birkhoff} obtained a  factorization theorem on matrix-valued functions which implies that any holomorphic vector bundle over the Riemann sphere is a direct sum of holomorphic line bundles. The non-homogeneous Riemann-Hilbert problem was considered in 1934 by Privalov \cite{Privalov}. The Riemann-Hilbert problem has a rich variety of applications in many fields; see in particular the monographs by Gakhov \cite{Gakhov}, Muskhelishvili \cite{Mus}, Vekua \cite{Vekua}, and Wegert \cite{Wegert}. 

We now return to Riemann's boundary value problem (\ref{eq:R}), but replacing the affine lines by closed Jordan curves $l_z\subset \c$ containing the origin in the bounded component of $\c\setminus l_z$ for every $z\in \t$. This {\em nonlinear Riemann-Hilbert problem} was solved by Forstneri\v c in 1988 \cite{F1988} who showed that the graphs of solutions fill the polynomial hull of the torus $T=\{(z,w)\in\c^2: z\in \t,\ w\in l_z\}$. (See also the papers \cite{AW,F1987,Snirelman}.) An extension to more general sets was given by Slodkowski \cite{Slod}. Berndtsson and Ransford \cite{BR} used these ideas to give another proof of Carleson's Corona Theorem for the disc. \v Cerne \cite{Cerne2004} found some solutions of the nonlinear Riemann-Hilbert problem on bordered Riemann surfaces.  

A closely related direction, which offers a more geometric point of view on Riemann's boundary value problem, is the existence and perturbation theory of analytic discs (and of bordered Riemann surfaces) with boundaries attached to a certain real submanifold in a complex manifold. This point of view was pioneered by Bishop \cite{Bishop1965} who studied local envelopes of holomorphy of real surfaces in $\c^2$ and, more generally, of real submanifolds near complex singularities, by introducing and solving the so-called {\em Bishop equation}. His work was continued and developed in several directions by Kenig and Webster \cite{KW}, Lempert \cite{Lempert1981}, Bedford and Gaveau \cite{BG}, Forstneri\v c \cite{F1987},  Tr\'epreau \cite{Trepreau}, Tumanov \cite{Tumanov1988,Tumanov1991,Tumanov1998}, Bedford and Klingenberg \cite{BK}, Globevnik \cite{Globevnik1994}, \v Cerne \cite{Cerne1995} and many others. This subject also received considerable attention on almost complex manifolds starting with Gromov's seminal work in 1985 \cite{Gromov1985} on the use of pseudoholomorphic curves in symplectic geometry. For these developments see e.g.\ the survey by Eliashberg \cite{Eliashberg} on the technique of filling by holomorphic discs, the collection by Audin and Lafontaine \cite{AudinL}, and the papers \cite{HLS,ST2008,ST2011}, among others. This line of work also plays an important role in modern topology, in particular in Floer homology. We are unable to present a comprehensive survey of this large body of results in a short space, so we apologize to the authors whose contributions are not mentioned in the above summary. 

We wish to emphasize that these are very difficult analytic problems whose exact solutions are typically extremely difficult or impossible to find. However, in many applications one only needs an {\em approximate solution}, and this is usually a much easier problem. The following version of the aproximate Riemann-Hilbert problem appears in many different constructions: of proper holomorphic maps, of bounded complete holomorphic curves, in formulas for extremal functions in pluripotential theory, etc. 

Let $X$ be a complex manifold (or a complex space). We are given a holomorphic map $f\colon \overline\d\to X$, also called an {\em analytic disc} in $X$, and for each point $z\in \t$ a holomorphic map $g_z\colon \overline\d\to X$ such that $g_z(0)=f(z)$ and the discs $g_z$ depend continuously on $z\in \t$. Set $T_z=g_z(\t) \subset X$ and $S_z=g_z(\overline \d)\subset X$ for $z\in\t$. Fix a distance function $\dist$ on $X$. Given numbers $0<r<1$ and $\epsilon>0$, the {\em approximate Riemann-Hilbert problem} asks for a holomorphic map $F\colon \overline\d\to X$ satisfying the following properties for some $r'\in [r,1)$:
\begin{itemize}
\item[\rm (a)] $\dist(F(z),T_z)<\epsilon$ for $z\in\t$,
\item[\rm (b)] $\dist(F(\rho z),S_z)<\epsilon$ for $z\in \t$ and $r'\le \rho\le 1$, and
\item[\rm (c)] $\dist(F(z),f(z))<\epsilon$ for $|z|\le r'$.
\end{itemize}

These conditions can be adapted to any bordered Riemann surface $\overline M$ instead of the disc $\overline \d$ as the domain of the maps $f$ and $F$. (Compare with the statement of Theorem \ref{th:RH}.  The domain of the maps $g_z$ is always the closed disc $\overline\d$.)

When $X=\c^n$, this approximate Riemann-Hilbert problem is easily solved as follows (see \cite{FG2001} or \cite{DF2012} for the details). Consider the map 
\[
	\t\times \overline\d \ni (z,w) \mapsto g_z(w) - f(z) \in\c^n
\]
which is continuous in $z$ and holomorphic in $w$. Note that this vanishes at $w=0$ for any $z\in \t$ since $g_z(0)=f(z)$. We can approximate it arbitrarily closely by a rational map 
\[
	G(z,w)= z^{-m} \sum_{j=1}^N A_j(z) w^j \in\c^n, 
\]
where the $A_j$'s are $\c^n$-valued holomorphic polynomials and $m\in \n$. Pick $k\in \n$
and set 
\begin{equation}\label{eq:F}
	F(z)= f(z)+ G(z,z^k)= f(z) + z^{k-m} \sum_{j=1}^N A_j(z) z^{k(j-1)}, \quad z\in \overline \d. 
\end{equation}
The pole at $z=0$ cancels if $k>m$, and one easily verifies that $F$ satisfies the properties (a)--(c) if the integer $k$ is chosen big enough.

It is not clear how to solve this problem in an arbitrary complex manifold $X$ and with the disc $\d$ (as the domain of $f$) replaced by a bordered Riemann surface $M$. However, in most applications (in particular, in those related to the present survey) it suffices to solve the problem on a small disc $\overline D \subset \overline M$ which intersects the boundary $bM$ in a compact arc $I\subset M$ around a given point $p\in bM$. Furthermore, replacing the disc $g_{p}\colon \overline\d\to X$ by its graph in $\overline \d \times X$ which has an open Stein neighborhood in $\c\times X$, we can reduce the local approximation problem over $\overline D$ to the standard case of discs in $\c^n$. 

To enable the gluing of a local solution (on $\overline D$) with the given map $f$ (to get a new map on $\overline M$) we actually solve the following modified Riemann-Hilbert problem. Pick a pair of smaller arcs $I_0,I_1\subset bM$ such that $p\in I_0\subset I_1\subset I$ and a cut-off function $\chi\colon bM\to [0,1]$ such that $\chi=1$ on $I_0$ and $\chi=0$ on $bM\setminus I_1$. Set $\tilde g_z(w)=g_z(\chi(z)\,w)$ for $z\in bM$ and $w\in \overline\d$. Note that $\tilde g_z$ agrees with $g_z$ for $z\in I_1$ and is the constant disc $w\to f(z)$ for any point $z\in bM\setminus I_1$. We define $\tilde g_z$ as the constant dics $f(z)$ for points $z\in bD \setminus I_1$. Let $\wt F:\overline D\to X$ be an approximate solution of the Riemann-Hilbert problem with the data $f|_{\overline D}$ and $\tilde g_z$, $z\in bD$. By choosing the integer $k$ in (\ref{eq:F}) big enough, $\wt F$ satisfies condition (a) for $z\in I_0$, it satisfies condition (b) for $z\in bD$, and it is uniformly close to $f$ on $\overline D \setminus U$ where $U\subset \overline M$ is any given neighborhood of the arc $I_1$. In particular, we can write $\overline M=A\cup B$ where $A,B\subset \overline M$ are closed smoothly bounded domains such that $A$ is the complement of a small neighborhood of the arc $I_1$, $B\subset \overline D$ contains a small neighborhood of $I_1$, we have the separation property $\overline{A\setminus B}\cap \overline{B\setminus A}=\emptyset$ (any such pair $(A,B)$ is called a {\em Cartain pair}), and $\wt F$ is uniformly close to $f$ on the domain $C=A\cap B$. 

At this point we wish to glue $f$ and $\wt F$. If $X=\c^n$, this is a {\em Cousin-I problem with bounds}: the difference $c=\wt F-f$ is holomorphic and small on $C=A\cap B$, and by solving the $\overline\partial$-equation with bounds on $\overline M$ it can be split as the difference $c=b-a$ where $a,b$ are holomorphic and small on $A$ and $B$, respectively. Then $\wt F-b=f-a$ on $C=A\cap B$ and hence the two sides amalgamate into a holomorphic map $F\colon\overline M\to X$ satisfying the appropriate analogues of the conditions (a)--(c). 

This simple method does not work in the absence of a linear structure on the manifold $X$. However, we can still reduce the problem to the linear case by the method of {\em gluing holomorphic sprays}. We give an outline and refer for the details (in this precise setting) to \cite{DF2012}. For the general method of gluing sprays see \cite[Chapter 5]{F2011}. 

We begin by embedding $f$ as the core map $f=f_0$ in a family of holomorphic maps $f_t\colon\overline M\to X$, depending holomorphically on a parameter $t\in U\subset\c^N$ in a neighborhood of $0\in\c^N$ for some $N\in\n$, such that the partial differential 
\[
	\frac{\partial f_t(z)}{\partial t} \bigg|_{t=0}: T_0\c^N\cong \c^N \to T_{f_0(z)} X
\]
is surjective for every point $z\in \overline M$. Such a family $\{f_t\}$ is called a {\em dominating (holomorphic) spray of maps}. Next we shrink $U$ around $0\in\c^N$ and solve the Riemann-Hilbert problem over $\overline D$ to get a holomorphic family $\wt F_t\colon \overline D \to X$ $(t\in U)$ approximating $f_t$ on $C=A\cap B$. If the approximation is close enough, there exists a holomorphic  map $C\times U\ni (x,t)\mapsto \gamma(x,t) \in\c^N$ close to the map $(x,t)\mapsto t$ such that 
\[
		f_t(x) =  \wt  F_{\gamma(x,t)}(x),\quad x\in C,\ t\in U.
\]
(The set $U$ shrinks again around $0$.) Now the difficult part is to split $\gamma$ in the form 
\[
	\gamma(x,\alpha(x,t))= \beta(x,t),\quad x\in C,\ t\in W, 
\]
where $W\subset U$ is a neighborhood of $0\in \c^N$ and $\alpha\colon A\times W\to \c^N$, $\beta\colon B\times W\to \c^N$ are holomorphic maps close to the map $(x,t)\mapsto t$. This is achieved by solving the Cousin-I problem with bounds and using the implicit function theorem in Banach spaces. Then
\[
    f_{\alpha(t,x)} = \wt F_{\beta(x,t)}(x),\quad x\in C,\ t\in W,
\]
so the two sides define a spray of maps $\overline M\to X$. By taking $t=0$ we get a map $F\colon \overline M\to X$ satisfying the desired properties provided that the approximations were sufficiently close.

Having explained the problem and the method of solving it, we now give a brief survey of applications of this techniqe to the construction of proper holomorphic maps. The use of this method in the construction of complete bounded holomorphic immersions of bordered Riemann surfaces is discussed in the following section.

Suppose that $M$ is a bordered Riemann surface and $f\colon \overline M\to \c^n$ is a holomorphic map where $n>1$. Assume that $0\notin f(bM)$, a condition which holds for a generic $f$. In order to push the boundary $f(bM)$ further towards infinity, we choose for every $z\in bM$ a unit vector $V(z) \in \c^n$ orthogonal to $f(z)$ and consider the linear disc $g_z\colon \overline \d \to\c^n$ defined by $g_z(w)=f(z)+wV(z)$. Furthermore, to localize the problem as explained above, we introduce a cut-off function $\chi\colon bM\to [0,1]$ with support on a small arc $I\subset bM$ and consider instead the maps 
\[
	g_z(w)= f(z) + w \delta \chi(z) V(z),\quad z\in bM,\ w\in \overline \d;
\]
here $\delta>0$ is a constant. If the arc $I$ is short enough (so that the image $f(z)$ for $z\in I$ does not vary very much), we can use a fixed vector $V$ orthogonal to the point $f(p)$ for some $p\in I$. A solution $F\colon \overline M\to \c^n$ of the Riemann-Hilbert problem with this data then approximates $f$ on most of $\overline M$, but on the arc $I_0\subset I$ where $\chi=1$ we have $F(z)\approx f(z)+\delta V(z)$, so $|F(z)|\ge |f(z)| + c\delta^2$ for some constant $c\in (0,1)$ close to $1$. We have thus increased the distance from the origin by a fixed amount on the arc $I_0\subset bM$. By systematically repeating this construction on finitely many arcs which cover $bM$ we increase the distance by a fixed amount on all of $bM$, while at the same time aproximating the map as closely as desired on a given compact subset of $M$. We can perform this operation inductively so that the resulting sequence $F_k\colon \overline M\to\c^n$ converges uniformly on compacta in $M$ to a proper holomorphic map $F=\lim_{k\to\infty} F_k\colon M\to \c^n$. 

This method easily adapts to the case when $\c^n$ is replaced by an arbitrary complex manifold $X$ of dimension $n>1$ which admits a smooth exhaustion function $\rho\colon X\to\r$ whose Levi form (which equals the complex Hessian in any system of local holomorphic coordinates on $X$) has at least two positive eigenvalues at every point of $X$. (Such a manifold is said to be $(n-1)$-complete; it is $(n-1)$-convex if the condition holds outside some compact subset of $X$. See Grauert \cite{Grauert} for the theory of $q$-convexity.) In this case one uses small holomorphic discs $g_z$ in $X$ $(z\in bM)$ lying in the zero locus of the Levi polynomial of the function $\rho$ at the point $f(z)\in X$ in the direction of a positive eigenvalue. Along this zero set the quadratic holomorphic term in the Taylor expansion of $\rho$ at $f(z)$ vanishes, so $\rho$ equals the Levi form plus terms of higher order, and therefore it increases quadratically along the image of $g_z$. (Taking the exhaustion function $\rho=|z_1|^2+\ldots+|z_n|^2$ on $\c^n$ gives  linear discs orthogonal to the given point of $\c^n$.) 

For applications of this method to the construction of proper holomorphic mappings from open Riemann surfaces see the papers \cite{Globevnik1988,Globevnik1989,FG1992,FG2001,Globevnik2000,DD2002,DD2004,DF2007}, listed chronologically. The paper \cite{DF2007} contains the most general results in this direction and also includes a survey of the previous work. This list does not include papers on {\em embedding} Riemann surfaces in $\c^2$ where different techniques are used; see e.g.\ \cite{FW1,FW2,GS,Majcen}. Some work has also been done on almost complex Stein manifolds of real dimension 4 \cite{CST}.

A similar technique is employed in the {\em Poletsky theory of discs} to obtain formulas expressing the envelopes of various disc functionals as pointwise minima over a family of analytic discs through a given point. This gives fairly explicit formulas for several extremal functions in pluripotential theory. The initial point was the remarkable discovery of Poletsky \cite{Poletsky1991,Poletsky1993} and Bu and Schachermayer \cite{BS} of the formula for computing the biggest plurisubharmonic minorant of an upper semicontinuous function as the envelope of the Poisson functional. This also gives Poletsky's characterization of polynomially convex hulls by sequences of analytic discs. For recent developments on this subject see the papers \cite{DF2012,DF2012-2} and the references therein.


\section{The Riemann-Hilbert problem for null holomorphic curves}
\label{sec:RH2}
Approximate solutions of Riemann-Hilbert boundary value problems, described in the previous section, have recently been used by the authors \cite{AF1} in the construction of proper {\em complete} holomorphic immersions of any bordered Riemann surface $M$ into the  ball of $\c^2$ (see Theorem \ref{th:CY-C2} above). Furthermore, in \cite{AF3} the Riemann-Hilbert technique was applied for the first time to the construction of complete bounded null curves (see Theorem \ref{th:CY-C3}) and of proper null curves in $\c^3$ with a bounded coordinate function (Theorem \ref{th:F3'}). We now describe the main ideas behind these developments. 

Let $f\colon \overline M\to \c^n$ be a holomorphic map such that $0\notin f(bM)$. Fix a point $p\in bM$. By pushing in the direction orthogonal to $f(p)\in \c^n$ for the amount $\delta>0$ (using the Riemann-Hilbert method) we increase the length of curves in $M$ ending near $p$ by approximately $\delta$, while the outer radius only increases by the order of $\delta^2$. Performing this construction recursively with a sequence $\delta_k>0$ such that $\sum_k \delta_k=\infty$ while $\sum_k\delta_k^2<\infty$ one therefore expects to get a bounded complete immersion (proper in a ball if one controls the procedure sufficiently carefully). However, there is a difficulty in controlling the distance estimate for divergent curves in $M$ on long boundary segments of $M$; undesired shortcurts may appear as is seen in the {\em sliding curtain model}. (Imagine that a curtain is held fixed at the lower end and is pulled horizontally at the upper end; this stretches each thread of the curtain, but the vertical distance between the edges remains the same.) 

To eliminate this problem, the authors in \cite{AF1} adjusted to this purpose the method of {\em exposing boundary points}, originally developed by Forstneri\v c and Wold \cite{FW1} in the construction of proper holomorphic embeddings of bordered Riemann surfaces into $\c^2$. One divides each of the boundary curves of $M$ into finitely many adjacent arcs $I_1,\ldots,I_m$ which are short enough so that the distance estimates for divergent curves in $M$ terminating on any of these arcs (and approaching the arc sufficiently `radially') can be controlled. Let $p_1,\ldots,p_m$ be the endpoints of this collection of arcs. Fix a small number $\delta>0$. In order to eliminate any shortcuts (after the deformation) which could be formed by curves in $M$ wandering almost horizontally in the vicinity of $bM$, we first attach to $f(\overline M)$ at the point $f(p_k)$ an embedded arc $\lambda_k\subset\c^n$ of length $>\delta$ which stays very close to $f(p_k)$. It is now possible to deform the image $f(M)$ so that the deformation is arbitrarily small away from the points $p_k$, while at $f(p_k)$ the image of $\overline M$ is stretched within a thin tube around $\lambda_k$ so that $f(p_k)$ goes to the other endpoint of $\lambda_k$. The effect of this deformation is that curves in $M$ which come sufficiently close to one of the points $p_k$ get elongated by at least $\delta$ (this happens in particular to curves ending on $bM$ near $p_k$), while the outer radius of the image almost does not increase. Now we complete the picture by applying the Riemann-Hilbert method on each of the arcs $I_k$ without destroying the effect of the first step. The cumulative effect of both deformations is that the lenght of any divergent curve in $M$ increases by at least $\delta>0$ while the outer radius only increases by $O(\delta^2)$. The proof is finished by a recursive procedure. The details are considerable (cf.\ \cite{AF1}). 

A similar construction can be done for null curves, except that the estimates become even more subtle (cf.\ \cite[\S 3]{AF3}). We now give a brief outline of this method.

\begin{proof}[Sketch of proof of Theorem \ref{th:RH}]
The special case $M=\d$ is done by an explicit calculation (cf.\ \cite[Lemma 3.1]{AF3}), using the so called {\em spinor representation} of the null quadric $\Agot$:
\[
	\pi\colon\c^2 \to \Agot,\quad 
	\pi(u,v)=\big(u^2-v^2,\imath (u^2+v^2),2uv\big).
\]
The restriction $\pi\colon \c^2\setminus\{0\} \to \Agot^*=\Agot\setminus\{0\}$ is an unbranched two-sheeted covering map, so the derivative $F'\colon \overline \d\to \Agot^*$ of any null disc $F\colon \overline \d\to \c^3$ lifts to a map $(u,v): \overline \d\to \c^2\setminus\{0\}$. One applies the Riemann-Hilbert problem to the map $(u,v)$ with a suitably chosen boundary data, projects the result by $\pi$ to the null quadric $\Agot^*$, and integrates to get a null disc $G \colon \overline \d\to \c^3$ satisfying properties (a), (b), and (c) in Theorem \ref{th:RH}, and also the following:
\begin{itemize}
\item[\rm (d)] $G$ is $\epsilon$-close to $F$ in the $\Cscr^1$ topology on $\overline{\d}\setminus U$, where $U\subset\overline\d$ is an arbitrarily small neighborhood of the arc $I\subset \t$.
\end{itemize}

To prove Theorem \ref{th:RH} in the general case we proceed as follows. Choose closed oriented Jordan curves $C_1,\ldots, C_m\subset M$ which determine a basis of the 1st homology group $H_1(M;\z)$. Pick a nowhere vanishing holomorphic $1$-form $\theta$ on $\overline M$ (such exists by the Oka-Grauert principle). Given a map $f\colon \overline M\to \c^3$, we denote by $\mathcal P_j(f)\in \c^{3}$ for $j=1,\ldots,m$ the period vector $\int_{C_j} f\theta$, and by $\mathcal P(f)\in \c^{3m}$ the period matrix with columns $\mathcal P_j(f)\in\c^3$. Observe that we have a bijective correspondence (up to constants)
\[
	\{F\colon \overline M\to\c^3\ \text{null\ curve}\} \longleftrightarrow
	\{f\colon \overline M\to \Agot^*\ \text{holomorphic},\ f\theta\  \text{exact}\}
\]
\[
		F(x)=F(p)+\int^x_p f\theta,\qquad dF=f\theta.
\]

Let $F\colon \overline M\to \c^3$ be a holomorphic null curve. Set $f=dF/\theta:\overline M\to \Agot^*$ and embed it as the core map $f=f_0$ in a dominating holomorphic spray of maps $f_t\colon \overline M\to \Agot^*$, where the parameter $t$ belongs to an open ball $B \subset\c^N$ around the origin.  Furthermore, we ensure that the associated period map $B\ni t\mapsto \mathcal P(f_t) \in \c^{3m}$ is submersive at $t=0$. 

Choose a small disc $D\subset M$ whose closure $\overline D\subset \overline M$ intersects $bM$ along a compact arc $J \subset bM$ which contains the given arc $I$ (containing the support of the function $\mu$ in Theorem \ref{th:RH}) in its relative interior. Fix a point $p\in D$ and apply the special case of Theorem \ref{th:RH} to the family of maps $F_t\colon \overline D\to \c^3$ given by $F_t(z)=F(p)+\int_{p}^z f_t\theta$ $(t\in B$) and the Riemann-Hilbert boundary data from Theorem \ref{th:RH}. (The integral is calculated along any path in $\overline D$ from $p$ to $z$; note that $F_0=F|_{\overline D}$.) This gives a new spray of null discs $\wt G_t\colon \overline D\to\c^3$ satisfying the properties (a)--(c) in Theorem \ref{th:RH} and also property (d) above (with $G$ replaced by $\wt G_t$ and $F$ replaced by $F_t$). The spray of maps $\tilde g_t= d\wt G_t/\theta: \overline M\to \Agot^*$ then approximates the initial spray $f_t\colon\overline D\to\Agot^*$ in the $\Cscr^1$ topology on $\overline D\setminus U$. If the approximations are close enough, we can glue these two sprays by the method described in \S\ref{sec:RH} to get a new holomorphic spray $g_t\colon\overline M\to \Agot^*$ for $t$ in a smaller ball $0\in B'\subset\c^N$. The spray $g_t$ approximates $f_t$ on $\overline M\setminus U$. Since we can choose the loops generating $H_1(M;\z)$ in this set, the period map $t\mapsto \mathcal P(g_t)$ approximates the map $t\mapsto\mathcal P(f_t)$. Since the latter map was chosen submersive at $t=0$, there exists a $t_0\in B'$ close to the origin such that the map $g=g_{t_0} \colon \overline M\to \Agot^*$ has vanishing periods. Hence $g$ integrates to a null curve $G(z)=F(p)+\int_p^z g\theta$ $(z\in\overline M)$. It can be verified that $G$  satisfies Theorem \ref{th:RH} provided that all aproximations were close enough. 
\end{proof}

The approximate solutions of the Riemann-Hilbert problem for null curves, established in \cite[\S3]{AF3}, have an additional feature which is the basis of the proof of Theorem \ref{th:F3}, and hence of the results in \S\ref{sec:SL2C}. When deforming a null disc $F\colon\overline \d \to\c^3$ in the direction of a null vector $V\in\Agot^*$ near a certain boundary point $p\in \t$ (see Theorem \ref{th:RH}), the component in the direction orthogonal to the 2-plane spanned by the vectors $F'(p) \in \Agot^*$ and $V$ changes arbitrarily little in the $\Cscr^0(\d)$ sense. To see how this is used, let us give

\begin{proof}[Sketch of proof of Theorem \ref{th:F3}] 
We outline the main idea in the simplest case when $M$ is the disc $\d$. Choose a couple of orthogonal null vectors in $\c^2\times\{0\}\subset\c^3$; for example, $V_1=(1,\imath,0)$ and $V_2=(1,-\imath,0)$. We begin with the linear null embedding $\overline \d \ni z\mapsto zV_1$. Using Theorem \ref{th:RH} with $M=\d$ we deform this embedding near the boundary in the direction of the vector $V_2$. Ignoring the nullity condition, one could simply take $\overline \d\ni z\mapsto zV_1+z^N V_2$ for a big integer $N$. To get a null curve, the third component must be involved, but it can be chosen arbitrarily $\Cscr^0$ small if $N$ is big enough. Next we deform the map from the previous step again in the direction of the vector $V_1$, changing  the third component only slighly in the $\Cscr^0(\d)$ norm. Repeating this alternating procedure gives a sequence of holomorphic null maps which converges to a proper null curve in $\c^3$ with a bounded third component. 

For a general bordered Riemann surface $M$ this construction is performed locally on small discs abutting the boundary $bM$, using also the method of exposed arcs in order to prevent any shorts. The local modifications are assembled together by the method of gluing sprays as described in the above sketch of proof of Theorem \ref{th:RH}.
\end{proof}

\subsection*{Acknowledgements}
A.\ Alarc\'{o}n is supported by Vicerrectorado de Pol\'{i}tica Cient\'{i}fica e Investigaci\'{o}n de la Universidad de Granada, and is partially supported by MCYT-FEDER grant MTM2011-22547 and Junta de Andaluc\'{i}a Grant P09-FQM-5088. 
F.\ Forstneri\v c is supported by the program P1-0291 and the grant J1-5432 from ARRS, Republic of Slovenia.


\vskip 0.2cm

\noindent Antonio Alarc\'{o}n

\noindent Departamento de Geometr\'{\i}a y Topolog\'{\i}a, Universidad de Granada, E-18071 Granada, Spain.

\noindent  e-mail: {\tt alarcon@ugr.es}

\vspace*{0.3cm}

\noindent Franc Forstneri\v c

\noindent Faculty of Mathematics and Physics, University of Ljubljana, and Institute
of Mathematics, Physics and Mechanics, Jadranska 19, 1000 Ljubljana, Slovenia.

\noindent e-mail: {\tt franc.forstneric@fmf.uni-lj.si}

\begin{thebibliography}{12}

\bibitem{Afer} Alarc\'on, A.; Fern\'andez, I.:
Complete minimal surfaces in $\r^3$ with a prescribed coordinate function. 
Differential Geom.\ Appl.\ \textbf{29} (2011), suppl.\ 1, S9--S15

\bibitem{AFL1} Alarc\'on, A.; Fern\'andez, I.; L\'opez, F.J.:
Complete minimal surfaces and harmonic functions. 
Comment.\ Math.\ Helv.\ \textbf{87} (2012)  891--904

\bibitem{AFL2} Alarc\'on, A.; Fern\'andez, I.; L\'opez, F.J.:
Harmonic mappings and conformal minimal immersions of Riemann surfaces into $\r^N$.
Calc.\ Var.\ Partial Differential Equations \textbf{47} (2013) 227--242

\bibitem{AFM} Alarc\'on, A.; Ferrer, L.; Mart\'in, F.:
Density theorems for complete minimal surfaces in $\r^3$. 
Geom.\ Funct.\ Anal.\ \textbf{18} (2008) 1--49

\bibitem{AF1}  Alarc\'{o}n, A.; Forstneri\v c, F.: 
Every bordered Riemann surface is a complete proper curve in a ball.
Math.\ Ann.\ \textbf{357} (2013) 1049--1070

\bibitem{AF2}  Alarc\'{o}n, A.; Forstneri\v c, F.: 
Null curves and directed immersions of open Riemann surfaces. 
Invent.\ Math., in press. \texttt{http://link.springer.com/article/10.1007/s00222-013-0478-8}

\bibitem{AF3}  Alarc\'{o}n, A.; Forstneri\v c, F.: 
The Calabi-Yau problem, null curves, and Bryant surfaces. 
Preprint 2013, \texttt{arXiv:1308.0903}

\bibitem{AL1} Alarc\'{o}n, A.; L\'{o}pez, F.J.:
Minimal surfaces in $\r^3$ properly projecting into $\r^2$. 
J.\ Diff.\ Geom. \textbf{90} (2012) 351--382

\bibitem{AL} Alarc\'{o}n, A.; L\'{o}pez, F.J.:
Null curves in $\c^3$ and Calabi-Yau conjectures. 
Math.\ Ann. \textbf{355} (2013) 429--455

\bibitem{AL-Israel} Alarc\'{o}n, A.; L\'{o}pez, F.J.:
Compact complete null curves in complex $3$-space.
Israel J.\ Math. \textbf{195} (2013) 97--122

\bibitem{AL3} Alarc\'{o}n, A.; L\'{o}pez, F.J.:
Proper holomorphic embeddings of Riemann surfaces with arbitrary topology into $\c^2$. 
J.\ Geom.\ Anal. \textbf{23} (2013) 1794--1805  

\bibitem{AL-conjugate} Alarc\'{o}n, A.; L\'{o}pez, F.J.:
Properness of associated minimal surfaces.
Trans.\ Amer.\ Math.\ Soc., in press. \texttt{arXiv:1203.0751}

\bibitem{AL-embedded} Alarc\'{o}n, A.; L\'{o}pez, F.J.:
A complete bounded embedded complex curve in $\c^2$. 
Preprint 2013, \texttt{arXiv:1305.2118}

\bibitem{AW} 
Alexander, H.; Wermer, J.: 
Polynomial hulls with convex fibers. Math.\ Ann. \textbf{271} (1985) 99--109

\bibitem{AudinL}
Lafontaine, J.; Audin, M. (eds.): 
Holomorphic curves in symplectic geometry. 
Progr. Math., 117, Birkh\"auser, Basel (1994)

\bibitem{BG}
Bedford, E.; Gaveau, B.:
Envelopes of holomorphy of certain 2-spheres in $\c^2$.
Amer.\ J.\ Math. \textbf{105} (1983) 975--1009 

\bibitem{BK}
Bedford, E.; Klingenberg, W.: 
On the envelope of holomorphy of a $2$-sphere in $\c^2$.
J.\ Amer.\ Math.\ Soc. \textbf{4} (1991) 623--646 

\bibitem{BN}
Bell, S.R.; Narasimhan, R.: 
Proper holomorphic mappings of complex spaces.
In \textit{Several complex variables, VI},
Encyclopaedia Math. Sci.\ \textbf{69}, 1--38,
Springer-Verlag, Berlin (1990)

\bibitem{BR} 
Berndtsson, B.; Ransford, T.J.: 
Analytic multifunctions, the $\overline\partial$-equation, and a proof of the corona theorem.
Pacific J.\ Math. \textbf{124} (1986) 57--72 

\bibitem{Birkhoff}
Birkhoff, G.D.: Singular points of ordinary linear differential equations.
Trans.\ Amer.\ Math.\ Soc. \textbf{10} (1909) 436--470

\bibitem{Bishop1965}
Bishop, E.:
Differentiable manifolds in complex Euclidean spaces.
Duke Math.\ J. \textbf{32} (1965) 1--21

\bibitem{Boli}
Bolibrukh, A.A.: Sufficient conditions for the positive solvability of the Riemann-Hilbert problem. Math.\  Notes \textbf{51} (1992) 110--117

\bibitem{Br} Bryant, R.: 
Surfaces of mean curvature one in hyperbolic space. 
Th\'eorie des vari\'et\'es minimales et applications (Palaiseau, 1983--1984).
Ast\'erisque \textbf{154-155} (1987), 12, 321--347, 353 (1988)

\bibitem{BS}
Bu, S.Q.; Schachermayer, W.:
Approximation of Jensen measures by image measures 
under holomorphic functions and applications.
Trans.\ Amer.\ Math.\ Soc. \textbf{331} (1992) 585--608

\bibitem{Calabi} Calabi, E.: 
Problems in differential geometry.
Ed.\ S.\ Kobayashi and J.\ Eells, Jr., Proceedings of the
United States-Japan Seminar in Differential Geometry, Kyoto, Japan, 1965. Nippon Hyoronsha Co.,
Ltd., Tokyo (1966) 170

\bibitem{CM} 
Colding, T.H.; Minicozzi II, W.P.: 
The Calabi-Yau conjectures for embedded surfaces. 
Ann.\ of Math. (2) \textbf{167}  (2008) 211--243

\bibitem{CHR} Collin, P.; Hauswirth, L.; Rosenberg, H.: 
The geometry of finite topology Bryant surfaces. 
Ann.\ of Math. (2) \textbf{153}  (2001) 623--659

\bibitem{CST}
Coupet, B.; Sukhov, A.; Tumanov, A.:
Proper J-holomorphic discs in Stein domains of dimension 2. 
Amer.\ J.\ Math. \textbf{131} (2009) 653--674 

\bibitem{Cerne1995}
\v Cerne, M.: Analytic discs attached to a generating CR-manifold. 
Ark.\ Mat. \textbf{33} (1995) 217--248 

\bibitem{Cerne2004}
\v Cerne, M.:
Nonlinear Riemann-Hilbert problem for bordered Riemann surfaces. 
Amer.\ J.\ Math.\ \textbf{126} (2004)  65--87

\bibitem{CF} 
\v Cerne, M.; Forstneri\v c, F.:
Embedding some bordered Riemann surfaces in the affine plane 
Math.\ Res.\ Lett. \textbf{9} (2002) 683--696 


\bibitem{DD2002}
Drinovec Drnov\v sek, B.:  
Proper holomorphic discs avoiding closed convex sets. 
Math.\ Z. \textbf{241} (2002) 593--596

\bibitem{DD2004}
Drinovec Drnov\v sek, B.: 
Proper discs in Stein manifolds avoiding complete pluripolar sets. 
Math.\ Res.\ Lett. \textbf{11} (2004) 575--581

\bibitem{DF2007} Drinovec Drnov\v sek, B.; Forstneri\v c, F.:
Holomorphic curves in complex spaces.
Duke Math.\ J. \textbf{139} (2007) 203--254

\bibitem{DF2012} Drinovec Drnov\v sek, B.; Forstneri\v c, F.:
The Poletsky-Rosay theorem on singular complex spaces.
Indiana Univ.\ Math.\ J. \textbf{61} (2012) 1707--1423 

\bibitem{DF2012-2} Drinovec Drnov\v sek, B.; Forstneri\v c, F.:
Disc functionals and Siciak-Zaharyuta extremal functions on singular varieties. 
Ann.\ Polon.\ Math. \textbf{106} (2012) 171--191 

\bibitem{Eliashberg}
Eliashberg, Y.:
Filling by holomorphic discs and its applications. In: 
Geometry of low-dimensional manifolds, 2 (Durham, 1989), 45--67,
London Math.\ Soc.\ Lecture Note Ser., 151, Cambridge Univ.\ Press, Cambridge (1990) 

\bibitem{Federer} Federer, H:
Geometric Measure Theory, Springer-Verlag, New York (1969)

\bibitem{FMM} Ferrer, L.; Mart\'{i}n, F.; Meeks, W.H.\ III: 
Existence of proper minimal surfaces of arbitrary topological type. 
Adv.\ Math. \textbf{231} (2012) 378--413

\bibitem{FMUY} Ferrer, L.; Mart\'{i}n, F.; Umehara, M.; Yamada, K.: 
A construction of a complete bounded null curve in $\c^3$.
Kodai Math.\ J., in press.

\bibitem{F1987-1} Forstneri\v c, F.: 
Analytic disks with boundaries in a maximal real submanifold of $\c^2$. 
Ann.\ Inst.\ Fourier (Grenoble) \textbf{37} (1987) 1--44

\bibitem{F1987} Forstneri\v c, F.:
On polynomially convex hulls with piecewise smooth boundaries. 
Math.\ Ann.\ \textbf{267} (1987) 97--104

\bibitem{F1988} Forstneri\v c, F.:
Polynomial hulls of sets fibered over the circle. 
Indiana Univ.\ Math.\ J. \textbf{37} (1988) 869--889

\bibitem{F2007} Forstneri\v c, F.: 
Manifolds of holomorphic mappings from strongly pseudoconvex domains. 
Asian J.\ Math. \textbf{11}  (2007) 113--126

\bibitem{F2011} Forstneri\v c, F.:  
Stein Manifolds and Holomorphic Mappings (The Homotopy Principle in Complex Analysis). 
Ergebnisse der Mathematik und ihrer Grenzgebiete, 3.\ Folge, 56. 
Springer-Verlag, Berlin-Heidelberg (2011)

\bibitem{FG1992} 
Forstneri\v c, F.; Globevnik, J.:
Discs in pseudoconvex domains.
Comment.\ Math.\ Helv.\ \textbf{67} (1992) 129--145

\bibitem{FG2001} 
Forstneri\v c, F.; Globevnik, J.:
Proper holomorphic discs in $\c^2$.
Math.\ Res.\ Lett. \textbf{8} (2001) 257--274

\bibitem{FW1} 
Forstneri\v c, F.; Wold, E.F.: 
Bordered Riemann surfaces in $\c^2$.  
J.\ Math.\ Pures Appl. \textbf{91} (2009) 100--114  

\bibitem{FW2} Forstneri\v c, F.; Wold, E.F.: 
Embeddings of infinitely connected planar domains into $\mathbb C^2$.
Anal.\ PDE \textbf{6} (2013) 499-–514

\bibitem{Gakhov} 
Gakhov, F.D.: Boundary value problems. 
Reprint of the 1966 English translation. 
Dover Publications, Inc., New York (1990) 

\bibitem{Globevnik1988}
Globevnik, J.: Boundary interpolation and proper holomorphic maps from the disc to the ball. 
Math.\ Z. \textbf{198} (1988) 143--150

\bibitem{Globevnik1989}
Globevnik, J.: Relative embeddings of discs into convex domains. 
Invent.\ Math. \textbf{98} (1989) 331--350

\bibitem{Globevnik1994}
Globevnik, J.: Perturbation by analytic discs along maximal real submanifolds of $\c^N$. 
Math.\ Z. \textbf{217} (1994) 287--316

\bibitem{Globevnik2000}
Globevnik, J.: Discs in Stein manifolds.
Indiana Univ.\ Math.\ J. \textbf{49} (2000) 553--574 

\bibitem{GS}
Globevnik, J.; Stens\o nes, B.:
Holomorphic embeddings of planar domains into $\c^2$.
Math.\ Ann. \textbf{303} (1995)  579--597


\bibitem{Grauert}
Grauert, H.: Theory of $q$-convexity and $q$-concavity.
In: Several complex variables, VII, 259--284, 
Encyclopaedia Math.\ Sci., 74, Springer-Verlag, Berlin (1994)

\bibitem{Gromov1985} Gromov, M.: 
Pseudoholomorphic curves in symplectic manifolds. 
Invent.\ Math. \textbf{82} (1985) 307--347

\bibitem{Hilbert}
Hilbert, D.:
Grundz\"uge einer allgemeinen Theorie der linearen Integralgleichungen.
Chelsea Publishing Company, New York, N.Y. (1953)

\bibitem{HLS}
Hofer, H.; Lizan, V.; Sikorav, J.-C.: 
On genericity for holomorphic curves in four-dimensional almost-complex manifolds. 
J.\ Geom.\ Anal. \textbf{7} (1997) 149--159 

\bibitem{HM} Hoffman, D.; Meeks III, W.H.:
The strong halfspace theorem for minimal surfaces.
Invent.\ Math.\ \textbf{101} (1990) 373--377 

\bibitem{Jones} Jones, P.W.: 
A complete bounded complex submanifold of $\c^3$. 
Proc.\ Amer.\ Math.\ Soc. \textbf{76} (1979) 305--306

\bibitem{JX} Jorge, L.P.; Xavier, F.:
A complete minimal surface in $\r^3$ between two parallel planes. 
Ann.\ of Math.\ (2) \textbf{112} (1980) 203--206

\bibitem{KW} 
Kenig, C.E.; Webster, S.M.: 
The local hull of holomorphy of a surface in the space of two complex variables. 
Invent. Math. \textbf{67} (1982) 1--21 

\bibitem{Kostov}
Kostov, V.P.: Fuchsian linear systems on $\c\mathbb P^1$ and the Riemann-Hilbert problem.
C.\ R.\ Acad.\ Sci.\ Paris S\'er.\ I Math. \textbf{315} (1992) 143--148

\bibitem{La} Lawson, B.:
Complete minimal surfaces in $S^3$.
Ann.\ of Math.\ (2) \textbf{92} (1970) 335--374

\bibitem{Lempert1981}
Lempert, L.:
La m\'etrique de Kobayashi et la repr\'esentation des domaines sur la boule. 
Bull.\ Soc.\ Math.\ France \textbf{109} (1981) 427–-474 

\bibitem{Lopez1} L\'opez, F.J.:
A nonorientable complete minimal surface in $\r^3$ between two parallel planes.
Proc.\ Amer.\ Math.\ Soc. \textbf{103} (1988) 913--917

\bibitem{Lopez2} L\'opez, F.J.:
Hyperbolic complete minimal surfaces with arbitrary topology.
Trans.\ Amer.\ Math.\ Soc. \textbf{350} (1998) 1977--1990

\bibitem{Majcen}
Majcen, I.:
Embedding certain infinitely connected subsets 
of bordered Riemann surfaces properly into $\c^2$.  
J.\ Geom.\ Anal. \textbf{19} (2009) 695--707 

\bibitem{MUY1} 
Mart\'{i}n, F.; Umehara, M.; Yamada, K.: 
Complete bounded null curves immersed in $\c^3$ and $SL(2,\c)$. 
Calc.\ Var.\ Partial Differential Equations \textbf{36} (2009) 119--139. 
Erratum: Complete bounded null curves immersed in $\c^3$ and $SL(2,\c)$.
Calc.\ Var.\ Partial Differential Equations \textbf{46} (2013) 439--440

\bibitem{MUY2} Mart\'{i}n, F.; Umehara, M.; Yamada, K.:
Complete bounded holomorphic curves immersed in $\c^2$ with arbitrary genus. 
Proc.\ Amer.\ Math.\ Soc. \textbf{137}  (2009) 3437--3450

\bibitem{Meeks} Meeks III, W.H:
Global problems in classical minimal surface theory. 
Global theory of minimal surfaces, 453--469, 
Clay Math. Proc., 2, Amer. Math. Soc., Providence (2005)

\bibitem{MPR} 
Meeks III, W.H.; P\'{e}rez, J.; Ros, A.: 
The embedded Calabi-Yau conjectures for finite genus. 

\texttt{http://www.ugr.es/$\sim$jperez/papers/papers.htm}

\bibitem{Morales} Morales, S.: 
On the existence of a proper minimal surface in $\r^3$ with a conformal type of disk. 
Geom.\ Funct.\ Anal.\ \textbf{13} (2003) 1281--1301

\bibitem{Mus} 
Muskhelishvili, N.I.: Singular integral equations. 
Corrected reprint of the 1953 English translation. Dover Publications, Inc., New York (1992)

\bibitem{Nad} Nadirashvili, N.:
Hadamard's and Calabi-Yau's conjectures on negatively curved and minimal surfaces.
Invent.\ Math. \textbf{126},  (1996) 457--465

\bibitem{Osserman}
Osserman, R.: A survey of minimal surfaces. Second ed. 
Dover Publications, Inc., New York (1986)

\bibitem{Plemelj}
Plemelj, J.:
Riemannsche Funktionenscharen mit gegebenen Monodromiegruppe.
Monatsh.\ Math.\ Phys. \textbf{19} (1908) 211--245

\bibitem{Plemelj2}
Plemelj, J.: Problems in the sense of Riemann and Klein. 
Edited and translated by J.R.M.\ Radok. 
Interscience Tracts in Pure and Applied Mathematics, 16. 
Interscience Publishers, John Wiley \& Sons Inc.,  New York-London-Sydney (1964) 

\bibitem{Poletsky1991}
Poletsky, E.A.: 
Plurisubharmonic functions as solutions of variational problems. 
In: Several complex variables and complex geometry, Part 1 
(Santa Cruz, CA, 1989), 163--171,
Proc.\ Sympos.\ Pure Math. \textbf{52}, Part 1, 
Amer.\ Math.\ Soc., Providence (1991) 

\bibitem{Poletsky1993}
Poletsky, E.A.: Holomorphic currents. 
Indiana Univ.\ Math.\ J. \textbf{42} (1993) 85--144

\bibitem{Privalov}
Privalov, I.I.:
On a boundary problem in analytic function theory. (Russian)
Mat.\ Sb. \textbf{41} (1934) 519--526 

\bibitem{Riemann}
Riemann, B.: Gesammelte mathematische Werke und wissenschaftlicher Nachlass. 
Dover Publications, Inc., New York (1953)

\bibitem{Ro} Rosenberg, H.: Bryant surfaces. 
In: The global theory of minimal surfaces in flat spaces (Martina Franca, 1999), 67-–111. 
Lecture Notes in Math., 1775, Springer, Berlin (2002)

\bibitem{RT} Rosenberg, H.; Toubiana, E.: 
A cylindrical type complete minimal surface in a slab of $\r^3$.
Bull.\ Sci.\ Math. (2) \textbf{111} (1987) 241--245

\bibitem{SY} Schoen. R.; Yau, S.T.: 
Lectures on Harmonic Maps.
Conference Proceedings and Lecture Notes in Geometry and Topology,
II. International Press, Cambridge, MA, 1997.

\bibitem{Slod}
Slodkowski, Z.: Polynomial hulls in $\c^2$ and quasicircles.
Ann.\ Scuola Norm.\ Sup.\ Pisa Cl.\ Sci.\ (4) \textbf{16} (1989) 367--391 (1990)

\bibitem{ST2008} Sukhov, A.; Tumanov, A.: 
Filling real hypersurfaces by pseudoholomorphic discs. 
J.\ Geom.\ Anal. \textbf{18} (2008) 632--649 

\bibitem{ST2011} Sukhov, A.; Tumanov, A.: 
Regularization of almost complex structures and gluing holomorphic discs to tori. 
Ann.\ Sc.\ Norm.\ Super.\ Pisa Cl.\ Sci. (5) \textbf{10} (2011) 389–-411 

\bibitem{Snirelman} \v Snirel'man, A.I.: 
The degree of a quasiruled mapping, and the nonlinear Hilbert problem. (Russian) 
Mat.\ Sb.\ (N.S.) \textbf{89 (131)} (1972) 366--389, 533 

\bibitem{Trepreau}
Tr\'epreau, J.-M.: 
Sur le prolongement holomorphe des fonctions CR d\'efinies sur une hypersurface 
r\'eelle de classe $\Cscr^2$ dans $\c^n$. Invent.\ Math. \textbf{83}  (1986) 583--592

\bibitem{Tumanov1988} 
Tumanov, A.E.: Extension of CR-functions into a wedge from a manifold of finite type. 
(Russian) Mat.\ Sb. (N.S.) \textbf{136 (178)} (1988) 128--139; translation in Math.\ USSR-Sb. \textbf{64} (1989) 129--140 

\bibitem{Tumanov1991} 
Tumanov, A.E.:
Extension of CR-functions into a wedge. (Russian) 
Mat.\ Sb.\ \textbf{181} (1990) 951--964; translation in Math. USSR-Sb. \textbf{70} (1991) 385--398

\bibitem{Tumanov1998} Tumanov, A.: 
Analytic discs and the extendibility of CR functions. 
Integral geometry, Radon transforms and complex analysis (Venice, 1996), 123--141,
Lecture Notes in Math., 1684, Springer, Berlin (1998)

\bibitem{UY} Umehara, M.; Yamada, K.: 
Complete surfaces of constant mean curvature $1$ in the hyperbolic $3$-space.
Ann.\ of Math. (2) \textbf{137}  (1993) 611--638

\bibitem{Vekua}
Vekua, N.P.: Systems of singular integral equations. 
P. Noordhoff, Ltd., Groningen (1967) 

\bibitem{Wegert}
Wegert, E.: 
Nonlinear boundary value problems for holomorphic functions and singular integral equations. 
Mathematical Research, 65. Akademie-Verlag, Berlin (1992)

\bibitem{Yang1} Yang, P.: 
Curvatures of complex submanifolds of $\c^n$. 
J.\ Diff.\ Geom. \textbf{12} (1977) 499--511

\bibitem{Yang2} Yang, P.:
Curvature of complex submanifolds of $\c^n$. 
Proc.\ Sympos.\ Pure Math., vol.\ 30, part 2, pp.\ 135--137.
Amer.\ Math.\ Soc., Providence (1977)

\bibitem{Yau1} Yau, S.T.: 
Problem section, Seminar on Differential Geometry.
Ann. of Math. Studies, vol.\ 102, pp.\ 669--706. Princeton University Press, Princeton (1982) 

\bibitem{Yau2} Yau, S.T.:
Review of geometry and analysis. 
Mathematics: frontiers and perspectives, 353--401. Amer. Math. Soc., Providence (2000)
 
\end{thebibliography}
\end{document}